\documentclass[11pt]{amsart}
\usepackage{amsmath,amssymb}
\newtheorem{theorem}{Theorem}
\newtheorem{proposition}[theorem]{Proposition}
\newtheorem{corollary}[theorem]{Corollary}
\newtheorem{lemma}[theorem]{Lemma}

\begin{document}

\title{Classification of manifolds with weakly $1/4$-pinched curvatures}
\author{Simon Brendle and Richard Schoen}
\address{Department of Mathematics \\
                 Stanford University \\
                 Stanford, CA 94305}
\thanks{The first author was partially supported by a Sloan Foundation Fellowship and by NSF grant DMS-0605223. The second author was partially supported by NSF grant DMS-0604960.}
\maketitle

\section{Introduction}

A classical theorem due to M.~Berger \cite{Berger2} and W.~Klingenberg \cite{Klingenberg1} states that a simply connected Riemannian manifold whose sectional curvatures all lie in the interval $[1,4]$ is either isometric to a symmetric space or homeomorphic to $S^n$ (see also \cite{Klingenberg2}, Theorems 2.8.7 and 2.8.10). In this paper, we provide a classification, up to diffeomorphism, of all Riemannian manifolds whose sectional curvatures are weakly $1/4$-pinched in a pointwise sense. Our main result is the following: 

\begin{theorem}
\label{weak.quarter.pinching}
Let $M$ be a compact Riemannian manifold of dimension $n \geq 4$. Suppose that $M$ has weakly $1/4$-pinched sectional curvatures in the sense that $0 \leq K(\pi_1) \leq 4 \, K(\pi_2)$ for all two-planes $\pi_1,\pi_2 \subset T_p M$. 
Moreover, we assume that $M$ is not locally symmetric. Then $M$ is diffeomorphic to a spherical space form.
\end{theorem}

In a previous paper \cite{Brendle-Schoen}, we proved that a compact Riemannian manifold with \textit{strictly} $1/4$-pinched sectional curvatures is diffeomorphic to a space form. Theorem \ref{weak.quarter.pinching} is a corollary of a more general theorem:

\begin{theorem} 
\label{main.theorem}
Let $M$ be a compact, locally irreducible Riemannian manifold of dimension $n \geq 4$. If $M \times \mathbb{R}^2$ has nonnegative isotropic curvature, then one of the following statements holds: 
\begin{itemize}
\item[(i)] $M$ is diffeomorphic to a spherical space form.
\item[(ii)] $n = 2m$ and the universal cover of $M$ is a K\"ahler manifold biholomorphic to $\mathbb{CP}^m$.
\item[(iii)] The universal cover of $M$ is isometric to a compact symmetric space.
\end{itemize}
\end{theorem}

By a theorem of R.~Hamilton \cite{Hamilton}, every compact, locally irreducible three-manifold with nonnegative Ricci curvature is diffeomorphic to a space form. Hence, Theorem \ref{main.theorem} also holds in dimensions two and three.

Theorem \ref{main.theorem} implies a structure theorem for compact Riemannian manifolds with the property that $M \times \mathbb{R}^2$ has nonnegative isotropic curvature. To explain this, suppose that $M$ is a compact Riemannian manifold such that $M \times \mathbb{R}^2$ has nonnegative isotropic curvature. By a theorem of Cheeger and Gromoll, the universal cover of $M$ is isometric to a product of the form $N \times \mathbb{R}^k$, where $N$ is a compact, simply-connected Riemannian manifold (cf. \cite{Cheeger-Gromoll} or \cite{Petersen}, p.~288). Moreover, $N$ is isometric to a product of the form $N_1 \times \hdots \times N_j$, where $N_1, \hdots, N_j$ are compact, simply connected, and irreducible (see \cite{Kobayashi-Nomizu1}, Chapter IV, Theorem 6.2). By Theorem \ref{main.theorem}, each of the factors $N_1, \hdots, N_j$ is either diffeomorphic to a sphere, or a K\"ahler manifold biholomorphic to complex projective space, or isometric to a compact symmetric space.

We now describe our strategy for handling the borderline case. Let $(M,g_0)$ be a compact Riemannian manifold of dimension $n \geq 4$, and let $g(t)$, $t \in [0,T)$, be the solution to the unnormalized Ricci flow with initial metric $g_0$. If $(M,g_0) \times \mathbb{R}^2$ has nonnegative isotropic curvature, then, by the results in \cite{Brendle-Schoen}, the product $(M,g(t)) \times \mathbb{R}^2$ has nonnegative isotropic curvature for all $t \in [0,T)$. If the manifold $(M,g(\tau))$ has general holonomy for some $\tau \in (0,T)$, we are able to use a strong maximum principle argument to show that the manifold $(M,g(\tau))$ satisfies the assumptions of Theorem 3 in our previous paper \cite{Brendle-Schoen}. This implies that $g_0$ will be deformed to a constant curvature metric by the normalized Ricci flow. In particular, $M$ is diffeomorphic to a spherical space form. 

We remark that there is a version of the strong maximum principle due to R. Hamilton \cite{Hamilton} for the curvature operator, but it does not seem to be sufficient for the present application. We believe the technique of this paper will have other applications to borderline situations which may be studied by Ricci flow methods.

In Section 2, we state a variant of the strong maximum principle for degenerate elliptic equations. This is a minor modification of a theorem of J.M.~Bony \cite{Bony}. 

In Section 3, we consider a family of metrics $g(t)$, $t \in [0,T)$, that have nonnegative isotropic curvature and evolve by the Ricci flow. For each time $t \in (0,T)$, we consider the set of all orthonormal four-frames with zero isotropic curvature. We show that this set is invariant under parallel transport. To that end, we view the isotropic curvature as a real-valued function on the frame bundle over $M$. Using results from \cite{Brendle-Schoen}, we show that this function satisfies a degenerate elliptic equation. This allows us to apply Bony's version of the strong maximum principle.

In Section 4, we complete the proof of Theorem \ref{main.theorem}. From this, Theorem \ref{weak.quarter.pinching} follows easily.

\section{A general maximum principle for degenerate elliptic equations}

Let $\Omega$ be an open subset of $\mathbb{R}^n$, and let $F \subset \Omega$ be relatively closed. As in \cite{Bony}, we say that a vector $\xi \in \mathbb{R}^n$ is tangential to $F$ at a point $x_1 \in F$ if $\langle x_1 - x_0,\xi \rangle = 0$ for all points $x_0 \in \mathbb{R}^n$ satisfying $d(x_0,F) = |x_1 - x_0|$. The following lemma is a slight modification of Theorem 2.1 in \cite{Bony}:

\begin{lemma}
\label{bony.lemma}
Let $\Omega$ be an open subset of $\mathbb{R}^n$, and let $F \subset \Omega$ be relatively closed. Assume that $X_1, \hdots, X_m$ are smooth vector fields on $\Omega$ that are tangential to $F$. Moreover, suppose that $\gamma: [0,1] \to \Omega$ is a smooth path such that $\gamma(0) \in F$ and $\gamma'(s) = \sum_{j=1}^m f_j(s) \, X_j(\gamma(s))$, where $f_1, \hdots, f_m: [0,1] \to \mathbb{R}$ are smooth functions. Then $\gamma(s) \in F$ for all $s \in [0,1]$. 
\end{lemma}

\textbf{Proof.} 
Choose a positive real number $\varepsilon$ such that $d(\gamma(s),\partial \Omega) \geq 2\varepsilon$ for all $s \in [0,1]$. We define a Lipschitz function $\rho: [0,1] \to \mathbb{R}$ by $\rho(s) = d(\gamma(s),F)^2$. We claim that there exists a positive constant $L$ such that \[\limsup_{h \searrow 0} \frac{1}{h} \, (\rho(s+h) - \rho(s)) \leq L \, \rho(s)\] for all $s \in [0,1)$ satisfying $\rho(s) \leq \varepsilon^2$. To see this, we fix a real number $s \in [0,1)$ such that $\rho(s) \leq \varepsilon^2$. There exists a point $x_1 \in F$ such that $d(\gamma(s),F) = |x_1 - \gamma(s)| \leq \varepsilon$. Since $X_1,\hdots,X_m$ are tangential to $F$ at $x_1$, we have $\langle x_1 - \gamma(s),X_j(x_1) \rangle = 0$ for $j = 1, \hdots, m$. This implies 
\begin{align*} 
&\limsup_{h \searrow 0} \frac{1}{h} \, (\rho(s+h) - \rho(s)) \\ 
&\leq \limsup_{h \searrow 0} \frac{1}{h} \, (|x_1 - \gamma(s+h)|^2 - |x_1 - \gamma(s)|^2) \\ 
&= -2 \, \sum_{j=1}^m f_j(s) \, \langle x_1 - \gamma(s),X_j(\gamma(s)) \rangle \\ 
&= 2 \, \sum_{j=1}^m f_j(s) \, \langle x_1 - \gamma(s),X_j(x_1) - X_j(\gamma(s)) \rangle \\ 
&\leq L \, |x_1 - \gamma(s)|^2 = L \, \rho(s) 
\end{align*} 
for some constant $L > 0$. Since $\rho(0) = 0$, we conclude that $\rho(s) = 0$ for all $s \in [0,1]$. 
This completes the proof. \\

\begin{proposition}
\label{bony.maximum.principle}
Let $\Omega$ be an open subset of $\mathbb{R}^n$, and let $X_1,\hdots,X_m$ be smooth vector fields on $\Omega$. Assume that $u: \Omega \to \mathbb{R}$ is a nonnegative smooth function satisfying 
\[\sum_{j=1}^m (D^2 u)(X_j,X_j) \leq -K \, \inf_{|\xi| \leq 1} (D^2 u)(\xi,\xi) + K \, |Du| + K \, u,\] 
where $K$ is a positive constant. Let $F = \{x \in \Omega: u(x) = 0\}$. Finally, let $\gamma: [0,1] \to \Omega$ be a smooth path such that $\gamma(0) \in F$ and $\gamma'(s) = \sum_{j=1}^m f_j(s) \, X_j(\gamma(s))$, where $f_1, \hdots, f_m: [0,1] \to \mathbb{R}$ are smooth functions. Then $\gamma(s) \in F$ for all $s \in [0,1]$. 
\end{proposition}

\textbf{Proof.} 
In view of Lemma \ref{bony.lemma}, it suffices to show that the vector fields $X_1, \hdots, X_m$ are tangential to $F$. In order to prove this, we adapt an argument due to J.M.~Bony (cf. \cite{Bony}, Proposition 3.1). Suppose that we are given two points $x_1 \in F$ and $x_0 \in \mathbb{R}^n$ such that $d(x_0,F) = |x_1 - x_0|$. We claim that $\langle x_1 - x_0,X_j(x_1) \rangle = 0$ for $j = 1, \hdots, m$. 

Without loss of generality, we may assume that $|x - x_0| > |x_1 - x_0|$ for all points $x \in F \setminus \{x_1\}$. (Otherwise, we replace $x_0$ by $\frac{1}{2} \, (x_0 + x_1)$.) Suppose that 
\[\sum_{j=1}^m \langle x_1 - x_0,X_j(x_1) \rangle^2 > 0.\] 
Then there exists a real number $\alpha > 0$ such that  
\[4\alpha^2 \, \sum_{j=1}^m \langle x_1 - x_0,X_j(x_1) \rangle^2 - 2\alpha \, \sum_{j=1}^m |X_j(x_1)|^2 > 2K\alpha + 2K\alpha \, |x_1 - x_0| + K.\] 
By continuity, there exists a bounded open set $U$ such that $x_1 \in U$, $\overline{U} \subset \Omega$, and 
\[4\alpha^2 \, \sum_{j=1}^m \langle x - x_0,X_j(x) \rangle^2 - 2\alpha \, \sum_{j=1}^m |X_j(x)|^2 > 2K\alpha + 2K\alpha \, |x - x_0| + K\] 
for all $x \in U$. As in \cite{Bony}, we define a function $v$ by 
\[v(x) = \exp(-\alpha |x - x_0|^2) - \exp(-\alpha |x_1 - x_0|^2).\] 
Moreover, we denote by $B$ the closed ball of radius $|x_1 - x_0|$ centered at $x_0$. By assumption, we have $u(x) > 0$ for all $x \in B \setminus \{x_1\}$. Since $\partial U \cap B$ is a compact subset of $B \setminus \{x_1\}$, there exists a real number $\lambda > 0$ such that $\lambda \, u(x) > v(x)$ for all $x \in \partial U \cap B$. Moreover, we have $\lambda \, u(x) \geq 0 > v(x)$ for all $x \in \partial U \setminus B$. Putting these facts together, we conclude that $\lambda \, u(x) > v(x)$ for all $x \in \partial U$. Pick a point $x_* \in \overline{U}$ such that $\lambda \, u(x_*) - v(x_*) \leq \lambda \, u(x) - v(x)$ for all $x \in U$. In particular, we have $\lambda \, u(x_*) - v(x_*) \leq \lambda \, u(x_1) - v(x_1) = 0$. Consequently, we have $x_* \in U$. This implies 
\begin{align*} 
&\sum_{j=1}^m (D^2 v)(X_j,X_j) \\ 
&\leq \lambda \, \sum_{j=1}^m (D^2 u)(X_j,X_j) \\ 
&\leq -K\lambda \, \inf_{|\xi| \leq 1} (D^2 u)(\xi,\xi) + K\lambda \, |Du| + K\lambda \, u \\ 
&\leq -K \, \inf_{|\xi| \leq 1} (D^2 v)(\xi,\xi) + K \, |Dv| + K \, v
\end{align*} 
at $x_*$. At the point $x_*$, we have 
\begin{align*} 
&\sum_{j=1}^m (D^2 v)(X_j,X_j) \\ 
&= \bigg [ 4\alpha^2 \, \sum_{j=1}^m \langle x_* - x_0,X_j(x_*) \rangle^2 - 2\alpha \, \sum_{j=1}^m |X_j(x_*)|^2 \bigg ] \, \exp(-\alpha |x_* - x_0|^2) 
\end{align*} 
and 
\[\inf_{|\xi| \leq 1} (D^2 v)(\xi,\xi) = -2\alpha \, \exp(-\alpha |x_* - x_0|^2).\] 
Moreover, we have $v(x_*) \leq \exp(-\alpha |x_* - x_0|^2)$. Putting these facts together, we obtain 
\[4\alpha^2 \, \sum_{j=1}^m \langle x_* - x_0,X_j(x_*) \rangle^2 - 2\alpha \, \sum_{j=1}^m |X_j(x_*)|^2 \leq 2K\alpha + 2K\alpha \, |x_* - x_0| + K.\] 
This contradicts our choice of $U$. \\

We point out that Proposition \ref{bony.maximum.principle} remains valid if $\Omega$ is a Riemannian manifold: 
to prove this, we subdivide the curve $\gamma$ into small segments, each of which is contained in a single coordinate chart. We then apply Proposition \ref{bony.maximum.principle} to each of these segments.

\section{Application to the Ricci flow}

In this section, we apply Bony's maximum principle to functions defined on the orthonormal frame bundle. Let $M$ be a compact manifold, and let $g(t)$, $t \in [0,T]$, be a family of metrics on $M$ evolving under Ricci flow. We consider the $(n+1)$-dimensional manifold $M \times (0,T)$. We denote by $E$ the vector bundle over $M \times (0,T)$ whose fiber over $(p,t) \in M \times (0,T)$ is given by $E_{(p,t)} = T_p M$. (In other words, $E$ is the pull-back of the tangent bundle $TM$ under the map $(p,t) \mapsto p$.) We define a bundle metric $h$ on $E$ by $\langle V,W \rangle_h = \langle V,W \rangle_{g(t)}$ for $V,W \in E_{(p,t)}$. We can extend the Riemannian connection to $M \times (0,T)$ by specifying the covariant time derivative. Given two sections $V,W$ of $E$, we define 
\begin{equation} 
\label{covariant.time.derivative}
\langle D_{\frac{\partial}{\partial t}} V,W \rangle_{g(t)} = \langle \frac{\partial}{\partial t} V,W \rangle_{g(t)} - \text{\rm Ric}_{g(t)}(V,W) 
\end{equation}
(cf. \cite{Chow-Knopf}). Of course we take $D_X V$ to be the Riemannian covariant derivative with respect to $g(t)$ for $X \in T_p M \subset T_{(p,t)}(M \times (0,T))$. This defines a connection on the vector bundle $E$ which is compatible with the bundle metric $h$.

We now define $P$ to be the orthonormal frame bundle of $E$ equipped with the natural right action of $O(n)$. Note that $P$ is a principal $O(n)$-bundle over $M \times (0,T)$. By definition, the fiber of $P$ over a point $(p,t) \in M \times (0,T)$ consists of all $n$-frames $\{e_1,\hdots,e_n\} \subset E_{(p,t)}$ that are orthonormal with respect to the bundle metric $h$. For each $t \in (0,T)$, we denote by $P_t$ the time $t$ slice of $P$. Clearly, $P_t$ is the orthonormal frame bundle of the Riemannian manifold $(M,g(t))$. 

For each $A \in O(n)$, we denote by $R_A$ the diffeomorphism on $P$ given by right translation by $A$. Given any $a \in \mathfrak{so}(n)$, we denote by $\sigma(a)$ the fundamental vertical vector field on $P$ whose value at $\underline{e} \in P$ is given by the tangent vector to the curve $s \mapsto R_{\exp(sa)}(\underline{e})$ at $s = 0$. The map $\sigma$ then defines a linear isomorphism from the Lie algebra $\mathfrak{so}(n)$ to the vertical subspace at $\underline{e}$. By transplanting the standard inner product on $\mathfrak{so}(n)$, we obtain a natural inner product on the vertical subspace at $\underline{e}$.

The connection $D$ defines a right invariant horizontal distribution on $P$. Hence, for each point $\underline{e} = \{e_1,\hdots,e_n\} \in P$, the tangent space $T_{\underline{e}} P$ can be written as a direct sum $\mathbb{H}_{\underline{e}} \oplus \mathbb{V}_{\underline{e}}$, where $\mathbb{H}_{\underline{e}}$ and $\mathbb{V}_{\underline{e}}$ denote the horizontal and vertical subspaces at $\underline{e}$, respectively. We next define a collection of smooth horizontal vector fields $\tilde{X}_1,\hdots,\tilde{X}_n,\tilde{Y}$ on $P$. For each $j = 1, \hdots, n$, the value of $\tilde{X}_j$ at a point $\underline{e} = \{e_1,\hdots,e_n\} \in P$ is given by the horizontal lift of the vector $e_j$. Moreover, the vector field $\tilde{Y}$ is defined as the horizontal lift of the vector field $\frac{\partial}{\partial t}$ on $M \times (0,T)$. Note that the vector fields $\tilde{X}_1, \hdots, \tilde{X}_n$ are tangential to $P_t$.

\begin{proposition} 
\label{maximum.principle.for.horizontal.laplacian}
Suppose that $u: P \to \mathbb{R}$ is a nonnegative smooth function which satisfies 
\begin{align*} 
&\tilde{Y}(u) - \sum_{j=1}^n \tilde{X}_j(\tilde{X}_j(u)) \\ 
&\geq K \, \inf_{\xi \in \mathbb{V}_{\underline{e}}, \, |\xi| \leq 1} (D^2 u)(\xi,\xi) - K \, \sup_{\xi \in \mathbb{V}_{\underline{e}}, \, |\xi| \leq 1} Du(\xi) - K \, u 
\end{align*} 
for some positive constant $K$. Let $F = \{u = 0\}$ be the zero set of $u$. Fix a real number $t \in (0,T)$, and let $\tilde{\gamma}: [0,1] \to P_t$ be a smooth horizontal curve satisfying $\tilde{\gamma}(0) \in F$. Then $\tilde{\gamma}(s) \in F$ for all $s \in [0,1]$.
\end{proposition}

\textbf{Proof.} 
Suppose that $\tilde{\gamma}: [0,1] \to P_t$ is a smooth horizontal curve satisfying $\tilde{\gamma}(0) \in F$. Then we can find smooth functions $f_1, \hdots, f_n: [0,1] \to \mathbb{R}$ such that $\tilde{\gamma}'(s) = \sum_{j=1}^n f_j(s) \, (\tilde{X}_j)_{\tilde{\gamma}(s)}$ for all $s \in [0,1]$. Hence, Proposition \ref{bony.maximum.principle} implies that $\tilde{\gamma}(s) \in F$ for all $s \in [0,1]$. \\

We now impose the additional condition that $(M,g(t))$ has nonnegative isotropic curvature for all $t \in [0,T]$. We define a nonnegative function $u: P \to \mathbb{R}$ by 
\begin{align*} 
u: \underline{e} = \{e_1, \hdots, e_n\} \mapsto \; &R(e_1,e_3,e_1,e_3) + R(e_1,e_4,e_1,e_4) \\ 
&+ R(e_2,e_3,e_2,e_3) + R(e_2,e_4,e_2,e_4) \\ &- 2 \, R(e_1,e_2,e_3,e_4), 
\end{align*} 
where $R$ denotes the curvature tensor of the evolving metric $g(t)$. The curvature tensor $R$ can be viewed as a section of the vector bundle $E^* \otimes E^* \otimes E^* \otimes E^*$. It follows from work of R.~Hamilton \cite{Hamilton} 
that 
\begin{equation} 
\label{evolution.of.curvature}
D_{\frac{\partial}{\partial t}} R = \Delta R + Q(R), 
\end{equation} 
where $D$ denotes the induced connection on the vector bundle $E^* \otimes E^* \otimes E^* \otimes E^*$, $\Delta$ is the Laplace operator with respect to the metric $g(t)$, and $Q(R)$ is a quadratic expression in the curvature tensor.

\begin{lemma} 
\label{diffusion.term}
At each point $\underline{e} = \{e_1,\hdots,e_n\} \in P$, we have 
\begin{align*} 
\tilde{Y}(u) - \sum_{j=1}^n \tilde{X}_j(\tilde{X}_j(u)) &= Q(R)(e_1,e_3,e_1,e_3) + Q(R)(e_1,e_4,e_1,e_4) \\ 
&+ Q(R)(e_2,e_3,e_2,e_3) + Q(R)(e_2,e_4,e_2,e_4) \\ 
&- 2 \, Q(R)(e_1,e_2,e_3,e_4). 
\end{align*} 
\end{lemma}

\textbf{Proof.} For each $j = 1, \hdots, n$, we have 
\begin{align*} 
\tilde{X}_j(\tilde{X}_j(u)) &= (D_{e_j,e_j}^2 R)(e_1,e_3,e_1,e_3) + (D_{e_j,e_j}^2 R)(e_1,e_4,e_1,e_4) \\ 
&+ (D_{e_j,e_j}^2 R)(e_2,e_3,e_2,e_3) + (D_{e_j,e_j}^2 R)(e_2,e_4,e_2,e_4) \\ &- 2 \, (D_{e_j,e_j}^2 R)(e_1,e_2,e_3,e_4). 
\end{align*} 
Summation over $j$ yields 
\begin{align*} 
\sum_{j=1}^n \tilde{X}_j(\tilde{X}_j(u)) &= (\Delta R)(e_1,e_3,e_1,e_3) + (\Delta R)(e_1,e_4,e_1,e_4) \\ 
&+ (\Delta R)(e_2,e_3,e_2,e_3) + (\Delta R)(e_2,e_4,e_2,e_4) \\ &- 2 \, (\Delta R)(e_1,e_2,e_3,e_4). 
\end{align*} 
Moreover, we have 
\begin{align*} 
\tilde{Y}(u) &= (D_{\frac{\partial}{\partial t}} R)(e_1,e_3,e_1,e_3) + (D_{\frac{\partial}{\partial t}} R)(e_1,e_4,e_1,e_4) \\ 
&+ (D_{\frac{\partial}{\partial t}} R)(e_2,e_3,e_2,e_3) + (D_{\frac{\partial}{\partial t}} R)(e_2,e_4,e_2,e_4) \\ &- 2 \, (D_{\frac{\partial}{\partial t}} R)(e_1,e_2,e_3,e_4). 
\end{align*} 
Hence, the assertion follows from (\ref{evolution.of.curvature}). \\

\begin{lemma} 
\label{reaction.term} 
At each point $\underline{e} = \{e_1, \hdots, e_n\} \in P$, we have 
\begin{align*} 
&Q(R)(e_1,e_3,e_1,e_3) + Q(R)(e_1,e_4,e_1,e_4) \\ 
&+ Q(R)(e_2,e_3,e_2,e_3) + Q(R)(e_2,e_4,e_2,e_4) \\ 
&- 2 \, Q(R)(e_1,e_2,e_3,e_4) \\ 
&\geq K \, \inf_{\xi \in \mathbb{V}_{\underline{e}}, \, |\xi| \leq 1} (D^2 u)(\xi,\xi) - K \, \sup_{\xi \in \mathbb{V}_{\underline{e}}, \, |\xi| \leq 1} Du(\xi) - K \, u, 
\end{align*} 
where $\mathbb{V}_{\underline{e}}$ denotes the vertical subspace at $\underline{e}$ and $K$ is a positive constant. 
\end{lemma}

\textbf{Proof.} 
We have 
\begin{align*} 
&Q(R)(e_1,e_3,e_1,e_3) + Q(R)(e_1,e_4,e_1,e_4) \\ 
&+ Q(R)(e_2,e_3,e_2,e_3) + Q(R)(e_2,e_4,e_2,e_4) - 2 \, Q(R)(e_1,e_2,e_3,e_4) \\ 
&= \sum_{p,q=1}^n (R_{13pq} - R_{24pq}) \, (R_{13pq} - R_{24pq}) \\ 
&+ \sum_{p,q=1}^n (R_{14pq} + R_{23pq}) \, (R_{14pq} + R_{23pq}) \\ 
&+ 2 \, I^{(1)} + 4 \, I^{(2)} + 2 \, I^{(3)}, 
\end{align*} 
where 
\begin{align*} 
I^{(1)} &= \sum_{p,q=1}^4 (R_{1p1q} + R_{2p2q}) \, (R_{3p3q} + R_{4p4q}) - \sum_{p,q=1}^4 R_{12pq} \, R_{34pq} \\ 
&- \sum_{p,q=1}^4 (R_{1p3q} + R_{2p4q}) \, (R_{3p1q} + R_{4p2q}) \\ &- \sum_{p,q=1}^4 (R_{1p4q} - R_{2p3q}) \, (R_{4p1q} - R_{3p2q}), 
\end{align*} 
\begin{align*} 
I^{(2)} &= \sum_{p=1}^4 \sum_{q=5}^n (R_{1p1q} + R_{2p2q}) \, (R_{3p3q} + R_{4p4q}) - \sum_{p=1}^4 \sum_{q=5}^n R_{12pq} \, R_{34pq} \\ 
&- \sum_{p=1}^4 \sum_{q=5}^n (R_{1p3q} + R_{2p4q}) \, (R_{3p1q} + R_{4p2q}) \\ 
&- \sum_{p=1}^4 \sum_{q=5}^n (R_{1p4q} - R_{2p3q}) \, (R_{4p1q} - R_{3p2q}), 
\end{align*} 
and 
\begin{align*} 
I^{(3)} &= \sum_{p,q=5}^n (R_{1p1q} + R_{2p2q}) \, (R_{3p3q} + R_{4p4q}) - \sum_{p,q=5}^n R_{12pq} \, R_{34pq} \\ 
&- \sum_{p,q=5}^n (R_{1p3q} + R_{2p4q}) \, (R_{3p1q} + R_{4p2q}) \\ 
&- \sum_{p,q=5}^n (R_{1p4q} - R_{2p3q}) \, (R_{4p1q} - R_{3p2q}). 
\end{align*}
By the results in Section 2 of \cite{Brendle-Schoen}, we can find positive constants $L_1,L_2,L_3$ such that 
\begin{align*} 
I^{(1)} &\geq -L_1 \, \sup_{\xi \in \mathbb{V}_{\underline{e}}, \, |\xi| \leq 1} Du(\xi) - L_1 \, u \\ 
I^{(2)} &\geq -L_2 \, \sup_{\xi \in \mathbb{V}_{\underline{e}}, \, |\xi| \leq 1} Du(\xi) \\ 
I^{(3)} &\geq L_3 \, \inf_{\xi \in \mathbb{V}_{\underline{e}}, \, |\xi| \leq 1} (D^2 u)(\xi,\xi) - L_3 \, u 
\end{align*}
for all $t \in (0,T)$. Putting these facts together, the assertion follows. \\

\begin{proposition} 
\label{isotropic.curvature}
Assume that $(M,g(t))$ has nonnegative isotropic curvature for all $t \in [0,T]$. Fix a real number $t \in (0,T)$. Then the set of all four-frames $\{e_1,e_2,e_3,e_4\}$ that are orthonormal with respect to $g(t)$ and satisfy 
\begin{align*} 
&R_{g(t)}(e_1,e_3,e_1,e_3) + R_{g(t)}(e_1,e_4,e_1,e_4) \\ 
&+ R_{g(t)}(e_2,e_3,e_2,e_3) + R_{g(t)}(e_2,e_4,e_2,e_4) \\ 
&- 2 \, R_{g(t)}(e_1,e_2,e_3,e_4) = 0 
\end{align*} 
is invariant under parallel transport.
\end{proposition}

\textbf{Proof.} 
Using Lemma \ref{diffusion.term} and Lemma \ref{reaction.term}, we obtain
\begin{align*} 
&\tilde{Y}(u) - \sum_{j=1}^n \tilde{X}_j(\tilde{X}_j(u)) \\ 
&\geq K \, \inf_{\xi \in \mathbb{V}_{\underline{e}}, \, |\xi| \leq 1} (D^2 u)(\xi,\xi) - K \, \sup_{\xi \in \mathbb{V}_{\underline{e}}, \, |\xi| \leq 1} Du(\xi) - K \, u. 
\end{align*}
Hence, the assertion follows from Proposition \ref{maximum.principle.for.horizontal.laplacian}. \\

We can draw a stronger conclusion if we assume that $(M,g(t)) \times \mathbb{R}^2$ has nonnegative isotropic curvature for all $t \in [0,T]$:

\begin{proposition} 
\label{product.with.R2}
Assume that $(M,g(t)) \times \mathbb{R}^2$ has nonnegative isotropic curvature for all $t \in [0,T]$. Fix real numbers $t \in (0,T)$ and $\lambda,\mu \in [-1,1]$. Then the set of all four-frames $\{e_1,e_2,e_3,e_4\}$ that are orthonormal with respect to $g(t)$ and satisfy 
\begin{align*} 
&R_{g(t)}(e_1,e_3,e_1,e_3) + \lambda^2 \, R_{g(t)}(e_1,e_4,e_1,e_4) \\ 
&+ \mu^2 \, R_{g(t)}(e_2,e_3,e_2,e_3) + \lambda^2 \mu^2 \, R_{g(t)}(e_2,e_4,e_2,e_4) \\ 
&- 2\lambda\mu \, R_{g(t)}(e_1,e_2,e_3,e_4) = 0 
\end{align*} 
is invariant under parallel transport.
\end{proposition}

\textbf{Proof.} 
We will apply Proposition \ref{isotropic.curvature} to the manifolds $(M,g(t)) \times S^1 \times S^1$, $t \in [0,T]$. 
Fix $\lambda,\mu \in [-1,1]$, and $t \in (0,T)$. Suppose that $\{e_1,e_2,e_3,e_4\} \subset T_p M$ is an orthonormal four-frame satisfying 
\begin{align*} 
&R(e_1,e_3,e_1,e_3) + \lambda^2 \, R(e_1,e_4,e_1,e_4) \\ 
&+ \mu^2 \, R_{g(t)}(e_2,e_3,e_2,e_3) + \lambda^2 \mu^2 \, R(e_2,e_4,e_2,e_4) \\ 
&- 2\lambda\mu \, R(e_1,e_2,e_3,e_4) = 0, 
\end{align*} 
where $R$ denotes the curvature tensor of $(M,g(t))$. We define an orthonormal four-frame $\{\hat{e}_1,\hat{e}_2,\hat{e}_3,\hat{e}_4\} \subset T_p M \times \mathbb{R}^2$ by 
\[\begin{array}{l@{\qquad\qquad}l} 
\hat{e}_1 = (e_1,0,0) & \hat{e}_2 = (\mu e_2,0,\sqrt{1-\mu^2}) \\ 
\hat{e}_3 = (e_3,0,0) & \hat{e}_4 = (\lambda e_4,\sqrt{1-\lambda^2},0). 
\end{array}\] 
The four-frame $\{\hat{e}_1,\hat{e}_2,\hat{e}_3,\hat{e}_4\}$ satisfies the relation 
\begin{align} 
\label{zero.isotropic.curvature}
&\hat{R}(\hat{e}_1,\hat{e}_3,\hat{e}_1,\hat{e}_3) + \hat{R}(\hat{e}_1,\hat{e}_4,\hat{e}_1,\hat{e}_4) \notag \\ 
&+ \hat{R}(\hat{e}_2,\hat{e}_3,\hat{e}_2,\hat{e}_3) + \hat{R}(\hat{e}_2,\hat{e}_4,\hat{e}_2,\hat{e}_4) \\ 
&- 2 \, \hat{R}(\hat{e}_1,\hat{e}_2,\hat{e}_3,\hat{e}_4) = 0, \notag
\end{align} 
where $\hat{R}$ denotes the curvature tensor of $(M,g(t)) \times \mathbb{R}^2$. It follows from Proposition \ref{isotropic.curvature} that the set of all orthonormal four-frames $\{\hat{e}_1,\hat{e}_2,\hat{e}_3,\hat{e}_4\}$ satisfying (\ref{zero.isotropic.curvature}) is invariant under parallel transport. This completes the proof.

\section{Proof of the main theorem}

Let $(M,g_0)$ be a compact Riemannian manifold of dimension $n \geq 4$ such that $(M,g_0) \times \mathbb{R}^2$ has nonnegative isotropic curvature. We denote by $g(t)$, $t \in [0,T)$, the solution to the Ricci flow with initial metric $g_0$. It follows from the results in \cite{Brendle-Schoen} that $(M,g(t)) \times \mathbb{R}^2$ has nonnegative isotropic curvature for all $t \in [0,T)$. In particular, $(M,g(t))$ has nonnegative sectional curvature for all $t \in [0,T)$.

\begin{proposition} 
\label{convergence}
Suppose that there exists a real number $\tau \in (0,T)$ such that $\text{\rm Hol}^0(M,g(\tau)) = SO(n)$. Then the normalized Ricci flow with initial metric $g_0$ exists for all time and converges to a constant curvature metric as $t \to \infty$.
\end{proposition} 

\textbf{Proof.} 
The assertion follows from Theorem 3 in \cite{Brendle-Schoen} if we can show that 
\begin{align*} 
&R_{g(\tau)}(e_1,e_3,e_1,e_3) + \lambda^2 \, R_{g(\tau)}(e_1,e_4,e_1,e_4) \\ 
&+ \mu^2 \, R_{g(\tau)}(e_2,e_3,e_2,e_3) + \lambda^2\mu^2 \, R_{g(\tau)}(e_2,e_4,e_2,e_4) \\ 
&- 2\lambda\mu \, R_{g(\tau)}(e_1,e_2,e_3,e_4) > 0 
\end{align*} 
for all orthonormal four-frames $\{e_1,e_2,e_3,e_4\}$ and all $\lambda,\mu \in [-1,1]$. In order to prove this, we fix a point $p \in M$ and real numbers $\lambda,\mu \in [-1,1]$. Suppose that $\{e_1,e_2,e_3,e_4\}$ is a four-frame in $T_p M$ which is orthonormal with respect to $g(\tau)$ and satisfies
\begin{align*} 
&R_{g(\tau)}(e_1,e_3,e_1,e_3) + \lambda^2 \, R_{g(\tau)}(e_1,e_4,e_1,e_4) \\ 
&+ \mu^2 \, R_{g(\tau)}(e_2,e_3,e_2,e_3) + \lambda^2 \mu^2 \, R_{g(\tau)}(e_2,e_4,e_2,e_4) \\ 
&- 2\lambda\mu \, R_{g(\tau)}(e_1,e_2,e_3,e_4) = 0. 
\end{align*} 
Since $\text{\rm Hol}^0(M,g(\tau)) = SO(n)$, the manifold $(M,g(\tau))$ is not flat. Hence, we can find a point $q \in M$ and an orthonormal two-frame $\{v_1,v_2\} \subset T_q M$ such that $R_{g(\tau)}(v_1,v_2,v_1,v_2) > 0$. Since $\text{\rm Hol}^0(M,g(\tau)) = SO(n)$, there exists a piecewise smooth path $\gamma: [0,1] \to M$ such that $\gamma(0) = p$, $\gamma(1) = q$, $v_1 = P_\gamma e_1$, and $v_2 = P_\gamma e_2$. (Here, $P_\gamma$ denotes parallel transport along $\gamma$ with respect to the metric $g(\tau)$.) Using Proposition \ref{product.with.R2}, we obtain 
\begin{align} 
\label{zero.isotropic.curvature.1} 
&R_{g(\tau)}(v_1,v_3,v_1,v_3) + \lambda^2 \, R_{g(\tau)}(v_1,v_4,v_1,v_4) \notag \\ 
&+ \mu^2 \, R_{g(\tau)}(v_2,v_3,v_2,v_3) + \lambda^2 \mu^2 \, R_{g(\tau)}(v_2,v_4,v_2,v_4) \\ 
&- 2\lambda\mu \, R_{g(\tau)}(v_1,v_2,v_3,v_4) = 0, \notag 
\end{align} 
where $v_3,v_4 \in T_q M$ are defined by $v_3 = P_\gamma e_3$ and $v_4 = P_\gamma e_4$. An analogous argument shows that 
\begin{align}
\label{zero.isotropic.curvature.2} 
&R_{g(\tau)}(v_1,v_2,v_1,v_2) + \lambda^2 \, R_{g(\tau)}(v_2,v_4,v_2,v_4) \notag \\ 
&+ \mu^2 \, R_{g(\tau)}(v_1,v_3,v_1,v_3) + \lambda^2 \mu^2 \, R_{g(\tau)}(v_3,v_4,v_3,v_4) \\ 
&- 2\lambda\mu \, R_{g(\tau)}(v_2,v_3,v_1,v_4) = 0 \notag 
\end{align} 
and 
\begin{align}
\label{zero.isotropic.curvature.3} 
&R_{g(\tau)}(v_2,v_3,v_2,v_3) + \lambda^2 \, R_{g(\tau)}(v_3,v_4,v_3,v_4) \notag \\ 
&+ \mu^2 \, R_{g(\tau)}(v_1,v_2,v_1,v_2) + \lambda^2 \mu^2 \, R_{g(\tau)}(v_1,v_4,v_1,v_4) \\ 
&- 2\lambda\mu \, R_{g(\tau)}(v_3,v_1,v_2,v_4) = 0. \notag 
\end{align} 
In the next step, we add equations (\ref{zero.isotropic.curvature.1}) -- (\ref{zero.isotropic.curvature.3}) and divide the result by $1 + \mu^2$. This yields 
\begin{align*} 
&\big [ R_{g(\tau)}(v_1,v_2,v_1,v_2) + R_{g(\tau)}(v_1,v_3,v_1,v_3) + R_{g(\tau)}(v_2,v_3,v_2,v_3) \big ] \\ 
&+ \lambda^2 \, \big [ R_{g(\tau)}(v_1,v_4,v_1,v_4) + R_{g(\tau)}(v_2,v_4,v_2,v_4) + R_{g(\tau)}(v_3,v_4,v_3,v_4) \big ] \\ &= 0. 
\end{align*} 
Since $(M,g(\tau))$ has nonnegative sectional curvature, it follows that \linebreak $R_{g(\tau)}(v_1,v_2,v_1,v_2) = 0$. This is a contradiction. \\

\begin{proposition} 
\label{product.with.R2.is.weakly.PIC}
Assume that $(M,g_0)$ is locally irreducible. Then one of the following statements holds: 
\begin{itemize}
\item[(i)] The normalized Ricci flow with initial metric $g_0$ exists for all time, and converges to a constant curvature metric as $t \to \infty$.
\item[(ii)] $n = 2m$ and the universal cover of $(M,g_0)$ is a K\"ahler manifold biholomorphic to $\mathbb{CP}^m$.
\item[(iii)] The universal cover of $(M,g_0)$ is isometric to a compact symmetric space.
\end{itemize}
\end{proposition}

\textbf{Proof.} 
By assumption, $(M,g_0)$ is compact, locally irreducible, and has nonnegative sectional curvature. Hence, a theorem of Cheeger and Gromoll (see \cite{Cheeger-Gromoll} or \cite{Petersen}, p.~288) implies that the universal cover of $M$ is compact. Suppose that $(M,g_0)$ is not locally symmetric. By continuity, there exists a real number $\delta \in (0,T)$ such that $(M,g(t))$ is locally irreducible and non-symmetric for all $t \in (0,\delta)$. According to Berger's holonomy theorem (see \cite{Joyce}, Theorem 3.4.1; see also \cite{Berger1},\cite{Simons}), there are three possibilities:

\textit{Case 1:} There exists a real number $t \in (0,\delta)$ such that $\text{\rm Hol}^0(M,g(t)) = SO(n)$. In this case, Proposition \ref{convergence} implies that the normalized Ricci flow with initial metric $g_0$ exists for all time and converges to a constant curvature metric as $t \to \infty$. 

\textit{Case 2:} $n = 2m$ and $\text{\rm Hol}^0(M,g(t)) = U(m)$ for all $t \in (0,\delta)$. In this case, the universal cover of $(M,g(t))$ is a K\"ahler manifold for all $t \in (0,\delta)$. Since $g(t) \to g_0$ in $C^\infty$, it follows that the universal cover of $(M,g_0)$ is a K\"ahler manifold. Moreover, the universal cover of $(M,g_0)$ is compact, irreducible, and has nonnegative sectional curvature. Hence, a theorem of N.~Mok \cite{Mok} implies that the universal cover of $(M,g_0)$ is either biholomorphic to $\mathbb{CP}^m$ or isometric to a Hermitian symmetric space. Since $(M,g_0)$ is not locally symmetric, the universal cover of $(M,g_0)$ must be biholomorphic to $\mathbb{CP}^m$.

\textit{Case 3:} $n = 4m \geq 8$ and there exists a real number $t \in (0,\delta)$ such that $\text{\rm Hol}^0(M,g(t)) = \text{\rm Sp}(m) \cdot \text{\rm Sp}(1)$. In this case, the universal cover of $(M,g(t))$ is a compact quaternionic-K\"ahler manifold with nonnegative sectional curvature. By a theorem of B.~Chow and D.~Yang, the universal cover of $(M,g(t))$ is isometric to a symmetric space (cf. \cite{Berger3},\cite{Chow-Yang}). This contradicts the fact that $(M,g(t))$ is non-symmetric. \\

\begin{corollary}
Assume that $(M,g_0)$ has weakly $1/4$-pinched sectional curvatures in the sense that $0 \leq K(\pi_1) \leq 4 \, K(\pi_2)$ for all two-planes $\pi_1,\pi_2 \subset T_p M$. Moreover, we assume that $(M,g_0)$ is not locally symmetric. Then the normalized Ricci flow with initial metric $g_0$ exists for all time, and converges to a constant curvature metric as $t \to \infty$.
\end{corollary}

\textbf{Proof.} 
By assumption, $(M,g_0)$ is not locally symmetric. In particular, $(M,g_0)$ is not flat. Since $(M,g_0)$ has weakly $1/4$-pinched sectional curvatures, it follows that $(M,g_0)$ is locally irreducible. By Theorem \ref{product.with.R2.is.weakly.PIC}, there are two possibilities: 

\textit{Case 1:} The normalized Ricci flow with initial metric $g_0$ exists for all time, and converges to a constant curvature metric as $t \to \infty$. In this case, we are done. 

\textit{Case 2:} $n = 2m$ and the universal cover of $(M,g_0)$ is a K\"ahler manifold. Since $(M,g_0)$ has weakly $1/4$-pinched sectional curvatures, the universal cover of $(M,g_0)$ is a K\"ahler manifold of constant holomorphic sectional curvature (cf. \cite{Kobayashi-Nomizu2}). Consequently, the universal cover of $(M,g_0)$ is isometric to complex projective space up to scaling. This contradicts the assumption that $(M,g_0)$ is not locally symmetric.

\end{document}